\def\qed{{\it Q.e.d.}}
\def\eps{\varepsilon}
\def\ergo{\Rightarrow}

\def\mlabel#1{\label{#1}}
\def\mref#1{(\ref{#1})}
\def\dfrac{\displaystyle\frac}
\def\ds{\displaystyle}
\def\nfloor#1{\left\lfloor #1 \right\rfloor}
\def\nceil#1{\left\lceil #1 \right\rceil}
\def\nround#1{\left\lfloor #1 \right\rceil}
\def\absolute#1{\left| #1 \right|}
\def\Z{\mathbb{Z}}
\def\N{\mathbb{N}}
\def\Q{\mathbb{Q}}
\def\R{\mathbb{R}}
\def\sign{\mathop{\hbox{sign}}}
\def\lg{\log_{10}}
\def\smallskip{\vskip 0.3em}

\documentclass[draft, 10pt]{article} 
\usepackage{amsfonts}

\newcounter{factno}
\def\thefactpun#1#2{\refstepcounter{factno}{\bf Fact~\thefactno\label{#1}#2\/}}
\def\thefact#1{\thefactpun{#1}{.}}
\def\fact#1{Fact~\ref{#1}}

\begin{document}

\title{A Rapidly Converging Machin-like Formula for $\pi$}
\author{Oleg S. Alferov}
\date{November 2023}
\maketitle

\begin{abstract}
We present a simple recurrent formula to generate the Machin-like expression
for calculating $\pi/4$. The method works for any denominator in the starting
term and always provides a finite decomposition.
We show that the terms in the Machin-like formula decrease so rapidly
that the Lehmer's measure can be made arbitrarily small only by selecting the first term.

We introduce the concept of the partial Machin-like formula.
While the growth of the integer numbers may quickly render the computer
implementation impractical, the same reason restricts the total contribution
of the high terms. If the required precision is known in advance, the
subset of the expression may be selected to satisfy it.

We also present the Python program to compute the terms of the Machin-like formula
(full and partial), and its Lehmer's measure.

\end{abstract}

\section{Background}

In the 1706, John Machin discovered the fact that the $\pi$ number can be expressed
as the finite sum of arctangents taken from certain fractions:
$$
\frac{\pi}{4} = 4 \arctan\frac{1}{5} - \arctan\frac{1}{239}.
$$

Since the Maclaurin series
$$
\arctan x = x - \frac{x^3}{3} + \frac{x^5}{5} - \frac{x^7}{7} + \dots
$$
converges fast if $x<1$ is small, Machin's discovery made it possible to calculate $\pi$
with high precision even before the advent of computers.

To honor his discovery, the identities in the form
\begin{equation}\mlabel{machin}
\frac{\pi}{4} = \sum_{k=0}^N m_k \arctan\frac{1}{q_k}, \quad
    \hbox{where $m_k \in \Z$ and $q_k \in \N$}
\end{equation}
are called the Machin-like formulas. Since then, many new identities have been
found. To estimate, in modern terms, the computational complexity of expressing $\pi$
by the formula~\mref{machin}, Lehmer introduced the following measure in 1938:
\begin{equation}\mlabel{lehmer}
\lambda = \sum_{k=0}^N \frac{1}{\lg q_k}.
\end{equation}

The smaller $\lambda$, the easier is the computation
of $\pi$.~\cite{abrar, source-pi, lehmer-cite}

\section{General facts and notations}

Let us remind readers of a few well-known facts and simple statements. For those who are
curious, either the formal proof is provided in the Supplement, or the reference is given
to the source where the fact can be found.

\subsection{Evolution of the arctangent argument}

Let $\tan\alpha =a/b$ and $\tan\beta = c/d$ be rational numbers.
Then $\tan(\alpha \pm \beta)$ is also the rational number,
immediately expressed via trigonometric identities.
Because of that, an arbitrary sum of arctangents from rational numbers is an arctangent from
the rational fraction again.~\cite{handbook}
\begin{equation}\mlabel{arctan-sum}
\begin{array}{l}
\tan(\alpha \pm \beta) = \dfrac{\tan\alpha \pm \tan\beta}{1 \mp \tan\alpha \tan\beta \strut},
    \quad\ergo \\[1em]
\arctan\dfrac{a}{b} \pm \arctan\dfrac{c \strut}{d} = \arctan\dfrac{ad \pm bc}{bd \mp ac}.
\end{array}
\end{equation}

The shortest form of the Machin-like formula~\mref{machin} provides the connection between
$\pi$ and arctangents:
\begin{equation}\mlabel{unity}
\frac{\pi}{4} = \arctan\frac{1}{1}.
\end{equation}

However,~\mref{unity} is the only single-term Machin-like formula.
\thefactpun{nosingle}:
$$
\frac{\pi}{4} \ne m \arctan \dfrac{1}{q_0}, \quad \hbox{$\forall$ $m$, $q_0 \in \N$, if $q_0>1$.}
$$

We are going to subtract arctangents of the selected fractions from the right-hand
side of~\mref{unity}, one by one.
Since $1 \in \Q$, the difference will be another arctangent of a rational number
on each step:
\begin{equation}\mlabel{remainder}
\dfrac{\pi}{4} = \arctan\dfrac{1}{1} =
    \left( \pm \arctan\dfrac{1}{q_0 \strut} + ... \pm \arctan\dfrac{1}{q_{n-1}} \right)
    + \arctan\dfrac{A_n}{B_n}.
\end{equation}

The sum in the brackets represents the Machin-like formula that we are building.
(If some of $q$-s repeat, we can group them into one term.)
The rightmost term outside brackets is the ``remainder''. If we add another term to the sum,
the ``remainder'' evolves:
\begin{equation}\mlabel{evolve}
\begin{array}{c}
\arctan\dfrac{A_n}{B_n} \mp \arctan\dfrac{1}{q_n}
    = \arctan\dfrac{A_{n+1}}{B_{n+1}}, \\[1.2em]
\left\{\begin{array}{l}
A_{n+1} = q_n A_n \mp B_n, \strut \\
B_{n+1} = q_n B_n \pm A_n,
\end{array}\right.
\end{array}
\end{equation}
and our expression for extending the Machin-like formula becomes
$$
\dfrac{\pi}{4} =
    \left( \pm \arctan\dfrac{1}{q_0 \strut} + ... \pm \arctan\dfrac{1}{q_{n-1}}
        \pm \arctan\dfrac{1}{q_n} \right)
    + \arctan\dfrac{A_{n+1}}{B_{n+1}}.
$$

We would also have to reduce the fraction by $\gcd(A_{n+1}, B_{n+1})$,
but it is irrelevant to our study.

The series starts with values $A_0=1$ and $B_0=1$.
On each step, we use $q_n$, extend the Machin-like expression with the next term,
and get the ``remainder'' corresponding to the pair $A_{n+1}$, $B_{n+1}$.
Our goal is to select the sequence
$\{ q_n \}$ in a way that the ``remainder'' becomes $0$ in the end, in other words,
$A_{N+1}=0$. After that, the iterations stop, and we get the new Machin-like formula.

\subsection{Nearest integer}

The floor function, the nearest integer less than or equal to the argument,
$y=\nfloor{x}$ is well defined.
In $\Q$, it ultimately comes from Euclid's division lemma~\cite{carmichael}.
(Let $x=a/b$, then $\exists$ $\alpha, \beta \in \Z$,
$0 \le \beta < b$, such as $a=\alpha b + \beta$, then $\nfloor{x} = \alpha$ and
$0 \le x - \nfloor{x} = \beta/b < 1$, \qed)
In $\R$, the floor can be introduced via approximations of $x$ in $\Q$.

The ceiling function, the nearest integer greater than or equal to the argument, can be defined
via the floor function: $\nceil{x} = -\nfloor{-x}$.
Properties of the floor and ceiling functions include:~\cite{knuth}
\begin{equation}\mlabel{floorceil}
\left\{\begin{array}{ll}
\nfloor{x} = x = \nceil{x}, & x \in \Z, \\[0.2em]
\nfloor{x} < x < \nceil{x} = \nfloor{x} + 1, \quad & x \notin \Z.
\end{array}\right.
\end{equation}

Interestingly, the notation for {\bf the nearest integer\/} is not so common. Let us use the
symbol $y=\nround{x}$ and the following intuitive definition.

{\bf Definition.\/} The nearest integer to $x \in \R$ is either $\nfloor{x}$, or $\nceil{x}$,
whichever is closer to $x$ (also, exact halves are rounded up):
\begin{equation}\mlabel{round-sel}
\nround{x} = \left\{
\begin{array}{ll}
\nfloor{x}, \quad & \absolute{\nfloor{x} - x} < \absolute{\nceil{x} - x}, \\[0.2em]
\nceil{x}, \quad & \hbox{otherwise.}
\end{array}
\right.
\end{equation}

\thefact{round}
$ \absolute{\nround{x} - x} \le 1/2$.

{\bf Corollary 1.\/}
If $x \in \R$ is given, and $y \in \Z$ is variable, the minimum of $\absolute{y-x}$
equals to $\absolute{\nround{x} - x}$.
(For proof, substitute $y = \nround{x} + v$.)

{\bf Corollary 2.\/}
Let $a, b \in R$ be given, $b>0$, and $m \in \Z$ is variable. Then $m = \nround{a/b}$
provides the minimum for expression $\absolute{a-mb}$.

\subsection{Miscellaneous}

For studying the partial Machin-like formulas, we need simplified means to estimate the higher terms in expression~\mref{machin}.

\thefact{atantaylor} Maclaurin series error for $x \in [0,1)$. See, for example,~\cite{thomas}.
\begin{equation}\mlabel{maclaurin2}
\arctan x = \sum_{k=0}^K (-1)^k \dfrac{k^{2k+1}}{2k+1} + \eps, \quad\hbox{and}\quad
    \absolute\eps \le \dfrac{x^{2K+3}}{2K+3}.
\end{equation}

{\bf Corollary.\/}
$x - x^3/3 < \arctan x < x$ for all $x \in (0,1)$.

\smallskip
\thefact{limits-exp}
Notable limit. See, for example,~\cite{calculus-free}.
$$
\lim_{x \to \infty} \left( 1 + \dfrac{a}{x} \right)^x = e^a
$$

\section{The idea}\label{idea}

Let us informally present the method to build a parametric Machin-like formula and postpone
the strict proof until the next section.
We select the starting $q_0$ and use identity~\mref{evolve} to subtract
the first term of~\mref{machin} until the difference is positive:
$$
\dfrac{\pi}{4} - m \arctan\dfrac{1}{q_0} = \arctan\dfrac{A_1}{B_1}, \quad\hbox{such as}\quad
0 < \dfrac{A_1}{B_1} < \dfrac{1}{q_0}.
$$

After the first term is established,
we extend the expression by adding the next terms.
We use the recurrent formula~\mref{evolve} again.
Let us keep $A_n \ge 0$ $\forall \, n$,
and require that the sequence $\{ A_n \}$ is strictly decreasing:
$$
0 \le A_{n+1} = q_n A_n - B_n < A_n.
$$

Given that $A_n > 0$, it is equivalent to
$$
\dfrac{B_n}{A_n} \le q_n < \dfrac{B_n}{A_n} + 1 \quad\ergo\quad 
q_n = \nceil{\dfrac{B_n}{A_n}}.
$$

Notice that $q_n > 0$ and $B_{n+1} > B_n > 0$, therefore $q_{n+1} > q_n$.
For this reason, expression~\mref{machin} has $m_0 = m$, and $m_k = 1$ for all $k>0$.

That is all. We calculate the next term based on the ``remainder'' until $A_{n+1} = 0$.
Since the sequence $\{ A_n \}$ of integers is limited from below and is strictly decreasing, it is finite.
Any starting $q_0 > 1$ generates a Machin-like identity in which terms
can be made arbitrarily small.

Unfortunately, this method may have slow convergence, and Lehmer's measure of
the resultant identity may be high.
We need stronger constraints to improve our identity.
 
\section{The series}

The series above uses all positive coefficients. Let us modify the approach and allow
negative terms in the series. We are now looking for the identity in the form
\begin{equation}\mlabel{machin-switch}
\dfrac{\pi}{4} = m \arctan \dfrac{1}{q_0} +
    \sum_{n=1}^N \delta_n \arctan \dfrac{1}{q_n},
    \quad\hbox{where $q_n > 0$, and $\delta_n=\pm 1$.}
\end{equation}

\subsection{Starting term}

Let $q_0 > 1$ be the first denominator.
Let us rewrite~\mref{evolve} for the case of subtracting the first arctangent term from $\pi/4$
multiple times:
\begin{equation}\mlabel{evolve-first-atan}
\dfrac{\pi}{4} - m \arctan\dfrac{1}{q_0} = \arctan\dfrac{a_m}{b_m}
\end{equation}
\vskip-1em
\begin{equation}\mlabel{evolve-first-int}
\left\{\begin{array}{l}
    a_0 = 1 \\
    b_0 = 1 \\
    a_{k+1} = q_0 a_k - b_k \\
    b_{k+1} = q_0 b_k + a_k \\
\end{array}\right.
\end{equation}

We need the absolute value of the ``remainder'' to be as small as possible.
According to the \fact{nosingle}, it will be nonzero.
With fixed $q_0$, because of \fact{round}, Corollary 2, it happens at the nearest integer of the fraction
\begin{equation}\mlabel{remainder-zero}
\gamma = \left. \dfrac{\pi}{4} \right/ \arctan\dfrac{1}{q_0}, \quad
m = \nround{ \gamma }
    \quad\ergo\quad
    0 < \absolute{ \arctan\dfrac{a_m}{b_m} } \le \dfrac{1}{2} \arctan\dfrac{1}{q_0}
\quad\hbox{ }
\end{equation}

For practical calculations, we cannot use real numbers or trigonometric functions.
Let us construct the method for computing $m$ using
the recurrent equations~\mref{evolve-first-int}.
Notice that the right-hand side of~\mref{evolve-first-atan} eventually becomes negative.
Specifically, the ``remainder'' is still positive at $m_{(-)} = \nfloor{\gamma}$, and negative
at $m_{(+)} = \nceil{\gamma}$. Proof. By the floor and ceiling properties~\mref{floorceil},
$$
0 < \gamma - m_{(-)} < 1 \quad\ergo\quad
0 < \dfrac{\pi}{4} - m_{(-)} \arctan\dfrac{1}{q_0} < \arctan\dfrac{1}{q_0}
$$
$$
m_{(+)} = m_{(-)} + 1 \quad\ergo\quad
\dfrac{\pi}{4} - m_{(+)} \arctan\dfrac{1}{q_0} < 0, \quad\hbox{\qed}
$$

Because the arctangent is monotonically increasing function,
$$
\arctan\dfrac{a_{m_{(-)}}}{b_{m_{(-)}}} > 0, \quad
\arctan\dfrac{a_{m_{(+)}}}{b_{m_{(+)}}} < 0 \quad\ergo\quad
\dfrac{a_{m_{(-)}}}{b_{m_{(-)}}} > 0, \quad
\dfrac{a_{m_{(+)}}}{b_{m_{(+)}}} < 0
$$

Notice from equations~\mref{evolve-first-int}, that when $a_{k+1}$ becomes
negative for the first time, all of $a_k$, $b_k$, and $b_{k+1}$ are still positive.
In other words, $m_{(-)}$ and $m_{(+)} = m_{(-)} + 1$ are defined by
the place where the sequence $\{a_k\}$ changes sign:
\begin{equation}\mlabel{cross-zero}
a_{m_{(-)}} > 0, \quad a_{m_{(+)}} < 0.
\end{equation}

Now, select between $m = m_{(-)}$ and $m = m_{(+)}$. The $(+)$ shall be
selected if the $(+)$ ``remainder'' is less than the $(-)$ one:
$$
\arctan\absolute{\dfrac{a_{m_{(-)}}}{b_{m_{(-)}}}} >
    \arctan\absolute{\dfrac{a_{m_{(+)}}}{b_{m_{(+)}}}}
\quad\ergo\quad m = m_{(+)},
$$
or, taking into account~\mref{cross-zero}, the same as
\begin{equation}\mlabel{m-zero}
m = \left\{
\begin{array}{ll}
m_{(+)}, \quad & a_{m_{(-)}} b_{m_{(+)}} + a_{m_{(+)}} b_{m_{(-)}} > 0, \\[0.2em]
m_{(-)}, & \hbox{otherwise.}
\end{array}
\right.
\end{equation}

\subsection{Reducing series}

We are building the Machin-like expression in form~\mref{machin-switch}, adding
terms one by one.
To simplify tracking the terms and ``remainders'' in~\mref{machin-switch}
with~\mref{evolve}, we introduce separate variable $\delta_n$ for the sign of the ``remainder''.
(The arctangent function is odd.)
Let us rewrite expression~\mref{remainder} for the ``remainder''
on each step:
\begin{equation}\mlabel{remainder-delta}
\dfrac{\pi}{4} = m \arctan\dfrac{1}{q_0}
    + \sum_{k=1}^{n-1} \delta_k \arctan\dfrac{1}{q_k}
    + \delta_n \arctan\dfrac{A_n}{B_n},
\end{equation}
where $\delta_n = \pm 1$ is selected to make positive $A_n \ge 0$ and $B_n > 0$.

The very first ``remainder''~\mref{evolve-first-atan} is
\begin{equation}\mlabel{remainder-first}
A_1 = \absolute{a_m}, \quad
B_1 = b_m, \quad
\delta_1 = \sign a_m.
\end{equation}

Let us rephrase~\mref{evolve} separating the fraction in the ``remainder'' from its sign.
The following equations also give us the expression for $\delta_n$ in~\mref{remainder-delta}:
\begin{equation}\mlabel{evolve-signed}
\begin{array}{c}
\mu_{n+1} \arctan\dfrac{A_{n+1}}{B_{n+1}} =
    \arctan\dfrac{A_n}{B_n} - \arctan\dfrac{1}{q_{n \strut}} \\
\left\{\begin{array}{l}
A_{n+1} = \absolute{q_n A_n - B_n \strut}, \\
\mu_{n+1} = \sign (q_n A_n - B_n), \\
B_{n+1} = q_n B_n + A_n, \\
\delta_{n+1} = \delta_n \mu_{n+1}.
\end{array}\right.
\end{array}
\end{equation}

Proof. Substitute the ``remainder'' in~\mref{remainder-delta}
according to~\mref{evolve-signed} and write
$$
\begin{array}{ll}
\dfrac{\pi}{4} = m \arctan\dfrac{1}{q_0}
    + \sum_{k=1}^{n-1} \delta_k \arctan\dfrac{1}{q_k}
    + \delta_n \arctan\dfrac{1}{q_n}
    + \delta_n \mu_{n+1} \arctan\dfrac{A_{n+1}}{B_{n+1}} \\[16pt]
\phantom{\dfrac{\pi}{4}} = m \arctan\dfrac{1}{q_0}
    + \sum_{k=1}^n \delta_k \arctan\dfrac{1}{q_k}
    + \delta_{n+1} \arctan\dfrac{A_{n+1}}{B_{n+1}}, \quad\hbox{\qed}
\end{array}
$$

If $A_n = 0$ in~\mref{remainder-delta} (or~\mref{evolve-signed}),
there is no ``remainder'', and we have the
identity. Let us assume that $A_n > 0$.
From the~\fact{round}, Corollary~2, minimum $A_{n+1}$ in~\mref{evolve-signed}
happens at
\begin{equation}\mlabel{q-next}
q_n = \nround{\dfrac{B_n}{A_n}}.
\end{equation}

Because of the \fact{round}, $A_{n+1}$ will be
\begin{equation}\mlabel{decay-a}
0 \le A_{n+1} = A_n \absolute{q_n - \dfrac{B_n}{A_n}} \le A_n / 2.
\end{equation}

Since~\mref{decay-a} implies that $A_{n+1} < A_n$, the sequence $\{ A_n \}$ is finite.
However, compared to all-positive Machin-like formula from the section~\ref{idea},
its signed version~\mref{machin-switch} converges faster, having as few as
$O(\ln A_1)$ terms. We will present a more meaningful estimation in the next
subsection.

\subsection{Characterization}

\subsubsection{The denominator $q_1$ of the second term in~\mref{machin-switch}}

Substitute~\mref{remainder-first} into~\mref{remainder-zero} and use~\mref{arctan-sum}
for double arctangent:
\begin{equation}\mlabel{compare-two-terms}
2\arctan \dfrac{A_1}{B_2} = \arctan \dfrac{2A_1/B_1}{1-\left(A_1/B_1\right)^{2\strut}} \le \arctan \dfrac{1}{q_0}
\end{equation}

Arctangent is a monotonically increasing function.
Taking into account the trivial corollary of~\mref{remainder-zero}
$$
\arctan \dfrac{A_1}{B_2} \le \dfrac{1}{2} \arctan \dfrac{1}{q_0} < \arctan \dfrac{1}{q_0},
$$
we can write $A_1/B_1 < 1/q_0 < 1$. Unequilibrium~\mref{compare-two-terms} transforms into
\begin{equation}\mlabel{first-remainder-frac}
\dfrac{2A_1}{B_1} \le \dfrac{1}{q_0} \left( 1 - \dfrac{A_1^2}{B_1^2} \right) < \dfrac{1}{q_0}
    \quad\ergo\quad \dfrac{B_1}{A_1} > 2 q_0.
\end{equation}

The second term is selected according to~\mref{q-next}. Use the \fact{round}:
$$
q_1 = \nround{ \dfrac{B_1}{A_1} } \ge \dfrac{B_1}{A_1} - 0.5 > 2 q_0 - 0.5,
    \quad q_0, q_1 \in \N \quad\ergo
$$
\begin{equation}\mlabel{second-term-double}
q_1 \ge 2 q_0
\end{equation}

\subsubsection{The growth of $\{ q_n \}$}

Let us estimate $q_{n+1}$ using its definition~\mref{q-next} and equations~\mref{evolve-signed},
\mref{decay-a}. Because of the nearest integer's properties,
$$
q_{n+1} \ge \dfrac{B_{n+1}}{A_{n+1}} - 0.5 \ge \dfrac{q_n B_n + A_n}{A_n / 2} - 0.5
    = 2 q_n \cdot \dfrac{B_n}{A_n} + 2 - 0.5.
$$
$$
q_n \le \dfrac{B_n}{A_n} + 0.5 \quad\ergo
$$
$$
\begin{array}{l}
q_{n+1} \ge 2 q_n ( q_n - 0.5 ) + 1.5 > 2 q_n^2 - q_n + 1 \\[0.5em]
\phantom{q_{n+1}} > 2 q_n^2 - 2 q_n + 1 = q_n^2 + (q_n - 1)^2 \ge q_n^2 \quad\ergo
\end{array}
$$
$$
q_{n+1} > q_n^2.
$$

Let us recalculate the estimation for $q_n$ via $q_1$.
{\bf Lemma.\/}
\begin{equation}\mlabel{estimate-qn1}
q_n > q_1 ^ {2^{n-1}} \quad \forall\, n \ge 2
\end{equation}

{\bf Proof\/} by induction. At $n=2$, observe:
$$
q_n = q_2 > q_1^2 = q_1^{2^{n-1}}.
$$

If the Lemma's statement holds at $n$, verify it at $n+1$:
$$
q_{n+1} > q_n^2 >\left( q_1^{2^{n-1}} \right)^2 = q_1^{2 \cdot 2^{n-1}} = q_1^{2^n},
    \quad\hbox{\qed}
$$

\subsubsection{The number of terms $N$ in~\mref{machin-switch}}

Since the overall length $N$ depends on the nominator $A_1$,
let us express it first.
From~\mref{evolve-first-int}, notice that
$$
\left\{ \begin{array}{l}
a_{k+1} < q_0 a_k \\
b_{k+1} > q_0 b_k
\end{array} \right.
    \quad\ergo\quad
\left\{ \begin{array}{l}
a_k < q_0^k \\
b_k > q_0^k
\end{array} \right.
    \quad\ergo\quad
b_k > a_k
$$

Use~\mref{evolve-first-int} again for $b_k$:
$$
b_{k+1} = q_0 b_k + a_k < (q_0 + 1) b_k
    \quad\ergo\quad b_k < (q_0 + 1)^k
$$

Taking into account the arctangent bounds, we estimate $m$:
$$
\begin{array}{l}
\arctan \dfrac{1}{q_0} = \dfrac{1}{q_0} \cdot \left( 1 + O \left(\dfrac{1}{q_0^2\strut} \right) \right),
    \quad\hbox{therefore, from~\mref{remainder-zero},} \\[1em]
m = \dfrac{\pi}{4} \left/ \arctan \dfrac{1}{q_0} \right. + O(1)
    = \dfrac{\pi}{4} \left/ \dfrac{1}{q_0} \right. \cdot
    \left( 1 + O \left( \dfrac{1}{q_0^2\strut} \right) \right) + O(1) \\[1em]
\phantom{m} = \dfrac{\pi q_0}{4} + O(1).
\end{array}
$$

Remember from~\mref{remainder-first}, \mref{first-remainder-frac},
that $B_1 = b_m$ and $A_1 / B_1 < 1 / (2 q_0)$.
Let us estimate $A_1$:
$$
\begin{array}{l}
A_1 = O \left( \dfrac{(q_0+1)^{\pi q_0 / 4 + O(1)}}{2 q_0} \right) \\[1em]
\phantom{A_1} = O \left( q_0^{\pi q_0 / 4} \cdot \dfrac{q_0^{O(1)}}{2q_0} \cdot
        \left( 1+\dfrac{1}{q_0} \right)^{\ds \pi q_0 / 4 + O(1/q_0)} \right) \\[1em]
\phantom{A_1} = O ( q_0^{\pi q_0 / 4} )
\end{array}
$$

Overall, taking into account the decay of $\{ A_n \}$~\mref{decay-a}, we can
estimate the total length of our Machin-like identity~\mref{machin-switch}:
\begin{equation}\mlabel{estimate-length}
N = O(\log_2 A_1) = O(\ln q_0^{\pi q_0 / 4}) = O(q_0 \ln q_0).
\end{equation}

\subsection{Lehmer's measure}

Let us summarize what we know so far~\mref{second-term-double},
\mref{estimate-qn1}, on the denominators of~\mref{machin-switch}
and then substitute that into~\mref{lehmer}.
\begin{equation}\mlabel{next-log}\vcenter{\hbox to 10em{\hss\vbox{
\begin{center}
\begin{tabular}{ l l }
$n=0:\quad$ & $\lg q_0$ \\
$n=1:$ & $\lg q_1 \ge \lg 2 + \lg q_0 > \lg q_0$ \\
$\dots$ \\
$n>1:$ & $\lg q_n > 2^{n-1} \lg q_1 > 2^{n-1} \lg q_0$ \\
\end{tabular}
\end{center}%
}\hss}}\end{equation}

Then,
\begin{equation}\mlabel{lehmer-switch}
\begin{array}{l}
\lambda = \dfrac{1}{\lg q_0} + \ds\sum_{n=1}^N \dfrac{1}{\lg q_n}
    < \dfrac{1}{\lg q_0} + \dfrac{1}{\lg q_0} \cdot \ds\sum_{n=1}^N 2^{-n+1} \\[1.5em]
\phantom{\lambda} < \dfrac{1}{\lg q_0} \cdot \left( 1 + 2\ds\sum_{n=1}^\infty 2^{-n} \right)
    = \dfrac{3}{\lg q_0}.
\end{array}
\end{equation}

So, the Lehmer's measure $\lambda$ can be made arbitrarily small by selecting
the appropriate large denominator $q_0$ of the first term.

\section{Partial series}

Let us consider the task of calculating $\pi$ with the precision~$\eps$ known in advance.
The terms in the Machin-like formula~\mref{machin-switch} decrease very fast.
According to~\mref{estimate-qn1} and~\mref{estimate-length}, the last term
can be very rougly estimated as $O(q_0^{-q_0^{\pi q_0/4}})$.
The chances are, that all the higher terms starting with a certain limit will be much less
than the required precision. Let us briefly discuss the strategy for selecting
the significant terms of~\mref{machin-switch} depending on~$\eps$.

For convenience, re-write~\mref{machin-switch} in form
\begin{equation}\mlabel{machin-x}
\dfrac{\pi}{4} = m X_0 + \ds\sum_{n=1}^N X_n, \quad
    \hbox{where $X_n = \delta_n \arctan\dfrac{1}{q_n}$ and $\absolute{\delta_n} = 1$.}
\end{equation}

Let us consider $0 < \eps_1 \ll 1$ such that for some $n$ fraction $1/q_n < \eps_1$.
Then
$$
\begin{array}{l}
\absolute{X_n} < 1/q_n < \eps_1, \\
\absolute{X_{n+1}} < \eps_1^2, \\
\absolute{X_{n+2}} < \eps_1^4 < \eps_1^3, \\
\dots \\
\absolute{X_{n+k}} < \eps_1^{k+1}, \\
\dots \\
\end{array}
$$

Therefore, the difference between the exact identity~\mref{machin-x} and the sum
of its beginning terms up to the $X_{n-1}$ will be
$$
\eps_{\mathrm{tail}} = \absolute{\ds\sum_{k=n}^N X_k} \le \ds\sum_{k=n}^N \absolute{X_k}
    < \eps_1 \ds\sum_{k=0}^\infty \eps_1^k = \dfrac{\eps_1}{1-\eps_1}.
$$

To be specific, let us introduce $\eps_2$ as function of $\eps_1$, $m$, and $n$:
$$
\eps_2 = \dfrac{\eps_1}{(1-\eps_1)(m+n-1)},
$$

The significant terms in~\mref{machin-x} are computed using the Maclaurin series for
arctangent~\mref{maclaurin2}. The error of calculating the arctangent function
depends on the last term $K$ in the arctangent decomposition.
Let us require that the error for each $i$-th term in~\mref{machin-x}
will be less than~$\eps_2$. The length $K_i$ of the Maclaurin decomposition
shall be large enough to satisfy the condition:
$$
\mathrm{error}_i \le \dfrac{1}{(2K_i+3) \, q_i^{2K_i+3\strut}} < \eps_2.
$$

The actual total error $\tilde\eps$ of our computation will be then
$$
\tilde\eps \le (m+n-1) \, \eps_2 + \eps_{\mathrm{tail}} = \dfrac{2 \eps_1}{1-\eps_1}.
$$

If we select $\eps_1$ based on the {\it required\/} precision $\eps$, we get the bound for
the actual error:
\begin{equation}\mlabel{partial-precision}
\eps_1 = \dfrac{\eps}{2+\eps} \quad\ergo\quad \tilde\eps \le \eps.
\end{equation}

Notice that all the variables involved in this discussion are functions of $q_0$ and $\eps$,
including the length $n$ of the partial Machin-like formula and the lengths of 
Maclaurin approximations for each accepted arctangent term.
Selection of $\eps_1$ and $\eps_2$ as shown above ensures that the maximum
calculation error will be not greater than~$\eps$, so we have constructed the practical
method to calculate~$\pi$.

Finally, Lehmer's measure for the incomplete series can be estimated using
inequality~\mref{next-log}. The last term of the partial Lehmer's sum is greater than
the sum of all the trailing Lehmer's terms.
(Assuming that the partial Machin-like formula has at least $2$ terms.)
Thus, doubling the last term in the partial Lehmer's sum
provides us with the upper bound for the actual Lehmer's measure.

\section{Numerical experiment}

Computational experiments were conducted for several starting $q_0$-s.
The implementation below uses equations~\mref{evolve-first-int},
\mref{m-zero}, \mref{evolve-signed}, and~\mref{next-log}.
Python~3.6 was used as the platform with long integer arithmetic support.
For partial Machin-like expressions we set the arbitrary limit of 1 million decimal digits.

\pagebreak
\begin{verbatim}
# Computation of recurrent Machin terms. Optionally: partial

import math;

def log10(x) :
    return math.log(x) / math.log(10);

def str_x_or_lg(nm, x) :   # def. lg(x) := log10(x)
    if x < 1e200 :
        return nm + " " + str(x);
    return "lg " + nm + " " + str(log10(x));

part = int(input("Partial ? (1=yes, 0=no) > "));
Q = int(input("Start Q > "));

A = 1;
B = 1;
QS = 1;
M = 0;
Q0 = Q;
QQ = [];
br = 0;

# first term
while A*Q >= B :
    print(str_x_or_lg("A", A), str_x_or_lg("B", B));

    M = M + 1;
    tmpA = A * Q - B;
    B = B * Q + A;
    A = tmpA;

# check other approximation for the first term
tmpA = A * Q - B;
tmpB = B * Q + A;
if A * tmpB + B * tmpA > 0 :
    QS = -1;
    M = M + 1;
    A = -tmpA;
    B = tmpB;
    print(str_x_or_lg("A", A), str_x_or_lg("B", B));

print("M", M, "\n---");

# higher terms
while A > 0 :
    print(str_x_or_lg("A", A), str_x_or_lg("B", B));

    Q = (B + A - 1) // A;

    tmpA = A * Q - B;
    tmpA1 = tmpA - A;

    if abs(tmpA) > abs(tmpA1) :
        Q = Q - 1;
        tmpA = tmpA1;

    B = B * Q + A;
    A = tmpA;

    QQ.append(QS * Q);  # trace the sign in q_k

    if A < 0 :
        A = -A;
        QS = -QS;

    if part > 0 and log10(Q) > 1000000 :
        print("break");
        br = 1;
        break;  # for partial series

# Lehmer measure
L = 1 / log10(Q0);
print("---\nM", M, str_x_or_lg("Q", Q0));
for q1 in QQ :
    L = L + 1 / log10(abs(q1));
    str_sign = ("(+)" if q1>0 else "(-)");
    print(str_sign, str_x_or_lg("Q", abs(q1)));

if br > 0 :  # remaining sum in the partial series is less 
    L = L + 1 / log10(abs(QQ[-1]));     # than the last term

# Pi sum for sanity check
S = M*math.atan(1/Q0);
for q1 in QQ :
    S = S + math.atan(1/q1);

say = ("(brk)\n---\nLehm <" if br>0 else "---\nLehm");
print(say, L, "\nPi", 4*S);
\end{verbatim}
\pagebreak

The program takes $q_0$ in the input and builds the corresponding Machin-like identity.
For the cases when the running variables become too large, there is an option
to terminate the calculation and make the formula partial. The program prints the starting
multiplier $m$, denominators $q_k$, the signs $\delta_k$, and the Lehmer's
measure $\lambda$ in the end. If a $q_k$ is too large for display (we selected the limit of $200$ digits),
its decimal logarithm is printed. For sanity check only, the computed $\pi$ is also printed.

A few outputs for various $q_0$ are listed below.

At $q_0 = 5$, the program prints the original Machin identity:
\begin{verbatim}
M 4 Q 5
(-) Q 239
---
Lehm 1.851127652316856
Pi 3.1415926535897936
\end{verbatim}

Let us provide the following identities to demonstrate how to interpret
the Python program output. The numbers can be copy-pasted into a software
capable of scientific computations such as Mathematica to validate the results of our study.

For $q_0 = 7$:
\begin{verbatim}
M 6 Q 7
(-) Q 15
(+) Q 1712
(-) Q 8886139
(+) Q 2526830931360443
---
Lehm 2.551666609279759 
Pi 3.1415926535897936
\end{verbatim}

The listing above corresponds to the expression
$$
\begin{array}{l}
6 \arctan\dfrac{1}{7} - \arctan\dfrac{1}{15} + \arctan\dfrac{1}{1712} \\[0.75em]
\phantom{0} - \arctan\dfrac{1}{8886139} + \arctan\dfrac{1}{2526830931360443}
\end{array}
$$

The code to enter to Wolfram Alpha\hfill\break
(link: {\tt https://www.wolframalpha.com/input\/})
\begin{verbatim}
Pi/4 == 6 ArcTan[1/7] - ArcTan[1/15] + ArcTan[1/1712]
    - ArcTan[1/8886139] + ArcTan[1/2526830931360443] 
\end{verbatim}

\pagebreak
For $q_0 = 8$.
Notice that while Python supports arbitrarily long integers, it {\it does not\/}
support arbitrary precision floating point numbers:
\begin{verbatim}
M 6 Q 8
(+) Q 25
(-) Q 1407
(+) Q 4150619
(+) Q 77950325308084
(+) Q 28355848339635153147414863515
(-) Q 2412162405181169014685016537064715579879917878585649329193
---
Lehm 2.4159383360928026 
Pi 3.141592653589793
\end{verbatim}
$$
\begin{array}{l}
6 \arctan\dfrac{1}{8}
+ \arctan\dfrac{1}{25}
- \arctan\dfrac{1}{1407}
+ \arctan\dfrac{1}{4150619} \\[0.75em]
+\, \arctan\dfrac{1}{77950325308084}
+ \arctan\dfrac{1}{28355848339635153147414863515} \\[0.75em]
-\, \arctan\dfrac{1}{2412162405181169014685016537064715579879917878585649329193}
\end{array}
$$

\begin{verbatim}

Pi/4 ==
6 ArcTan[1/8] + ArcTan[1/25] - ArcTan[1/1407] + ArcTan[1/4150619]
    + ArcTan[1/77950325308084]
    + ArcTan[1/28355848339635153147414863515]
    - ArcTan[1/24121624051811690146850165370
               64715579879917878585649329193]
\end{verbatim}

For $q_0=9$:
\begin{verbatim}
M 7 Q 9
(+) Q 93
(+) Q 22055
(+) Q 5085558009
(+) Q 767266041127734416424
(+) Q 1766091533603478722982708121680411788426907
---
Lehm 1.9607629078499424 
Pi 3.1415926535897927

Pi/4 == 7 ArcTan[1/9] + ArcTan[1/93] + ArcTan[1/22055]
    + ArcTan[1/5085558009] + ArcTan[1/767266041127734416424]
    + ArcTan[1/1766091533603478722982708121680411788426907]
\end{verbatim}

\pagebreak
For $q_0=10$:
\begin{verbatim}
M 8 Q 10
(-) Q 84
(-) Q 21342
(-) Q 991268848
(-) Q 193018008592515208050
(-) Q 197967899896401851763240424238758988350338
(-) Q 117573868168175352930277752844194126767991
      915008537018836932014293678271636885792397
---
Lehm 1.9473700443296986 
Pi 3.141592653589794

Pi/4 == 8 ArcTan[1/10] - ArcTan[1/84] - ArcTan[1/21342]
    - ArcTan[1/991268848] - ArcTan[1/193018008592515208050]
    - ArcTan[1/197967899896401851763240424238758988350338]
    - ArcTan[1/117573868168175352930277752844194126767991
               915008537018836932014293678271636885792397]
\end{verbatim}

Related to the famous approximation of $\pi \approx 22/7$ (see, for example,~\cite{source-pi}),
expression for $q_0 = 28$. For the sake of performance, numbers longer than $200$ decimal digits
are represented as their decimal logarithms. While, strictly speaking, the listing below is no longer
an identity, one can reproduce this test and output exact integer numbers for this Machin-like
formula.
Notice that $\log_{10} q_n$ almost exactly doubles on each step.
\begin{verbatim}
M 22 Q 28
(+) Q 56547
(+) Q 20747394343
(+) Q 1112172624652580034840
(-) Q 16659543628852678157467292276729792021493732
(+) Q 19351587917741573692739048650182250035782554284801229
      80428023197249578624178441690588894
(+) Q 14718492206740001931852838656976183022784010091410392
      42953147036168205460675285916208006732990521412670908
      69513168086930986444104325857945434713227531064709901
      94861973862674124
(-) lg Q 350.7305238264204
(-) lg Q 702.0893561664352
(+) lg Q 1404.5900031211877
(+) lg Q 2809.9358190450657
(-) lg Q 5620.463702225073
(+) lg Q 11241.25183905937
(-) lg Q 22484.181013176003
(-) lg Q 44968.75144493231
(-) lg Q 89937.82819599868
(-) lg Q 179876.09422636102
(+) lg Q 359752.6872249542
(+) lg Q 719508.3122952792
(-) lg Q 1439017.5723335177
(-) lg Q 2878035.9207072803
(-) lg Q 5756072.228487223
(-) lg Q 11512146.246898009
---
Lehm 1.091872372535026 
Pi 3.141592653589793
\end{verbatim}

To illustrate how Lehmer's measure decreases with the first term decrease, the following is
the example with $q_0 = 100000$. Partial mode is on.
Notice that~\mref{lehmer-switch} provides a close but
somewhat larger bound $\lambda < 0.6$.

\begin{verbatim}
M 78540 Q 100000
(-) Q 544491
(+) Q 783664894308
(+) Q 1303088915612811138696591
(+) Q 7636018810382840305552700218709810164960367081459
(+) Q 3617852367571966985355152435997587992571879865162469
      49685961215997793793692115231994619217388993130
(-) Q 1263548033106645763782664751160149068356817307489653
      3051095395918146583179586035763817403459556538268035
      7855957343101646766103100085828143285482786260424809
      4974093639535779278903243063902741055991584
(+) lg Q 396.72088863680796
(-) lg Q 793.8533155269043
(-) lg Q 1588.4525139699301
(+) lg Q 3177.6648529734907
(-) lg Q 6356.3408844093965
(+) lg Q 12713.353728781887
(-) lg Q 25427.097270768576
(-) lg Q 50855.26878265154
(+) lg Q 101710.88560659182
(+) lg Q 203422.37580891087
(-) lg Q 406845.1926757058
(-) lg Q 813691.1891000423
(-) lg Q 1627383.4447412174
(brk)
---
Lehm < 0.5405713556036384 
Pi 3.141592653589794
\end{verbatim}

\section{Supplement}

\noindent{\bf \fact{nosingle}.\/} Let us provide the formal proof that
the single-term identity~\mref{machin} is only possible if $q_0=1$.

\noindent{\bf Proof\/} by contradiction. Let $q_0>1$, and the single-term identity is possible.
Since $q_0>0$, we have $m>0$:
$$
\arctan \dfrac{1}{1} - \underbrace{\arctan \dfrac{1}{q_0} - \dots - \arctan \dfrac{1}{q_0}}_{\hbox{$m$ times}} = 0.
$$

The recurrent formula~\mref{evolve} defines the sequences $\{ A_n \}$ and $\{ B_n \}$ that
describe the result of the subtractions.
Because of their definition, $A_n$ and $B_n$ are polynomials of variable $q_0$
with whole coefficients. By our hypothesis, $q_0$ is the root of $A_m(q_0)=0$.

Let us trace how the constant term evolves in $A_n$ and $B_n$.
We have
$$
\begin{array}{l}
\left\{\begin{array}{l}
A_n(q_0) = a_n^n q_0^n + \dots + a_n^0, \\
B_n(q_0) = b_n^n q_0^n + \dots + b_n^0,
\end{array}\right. \\[1em]
\left\{\begin{array}{l}
A_0=1, \\
B_0=1,
\end{array}\right.
\quad\ergo\quad
\left\{\begin{array}{l}
a_0^0 = 1, \\
b_0^0 = 1,
\end{array}\right. \\[1em]
\left\{\begin{array}{l}
A_{n+1}(q_0) = q_0 A_n(q_0) - B_n(q_0) = a_n^n q_0^{n+1} + \dots - b_n^0, \\
B_{n+1}(q_0) = q_0 B_n(q_0) + A_n(q_0) = b_n^n q_0^{n+1} + \dots + a_n^0,
\end{array}\right.  \quad\ergo \\[1em]
\left\{\begin{array}{l}
a_{n+1}^0 = - b_n^0, \\
b_{n+1}^0 = a_n^0.
\end{array}\right.
\end{array}
$$

The constant coefficients change according to the table:
\begin{center}
\begin{tabular}{ c | c | c }
$n$ & $a_n^0$ & $b_n^0$ \\
\hline
0 & 1 & 1 \\
1 & -1 & 1 \\
2 & -1 & -1 \\
3 & 1 & -1 \\
4 & 1 & 1\hbox to 0pt{$\quad$loops to $n=0$\hss} \\
5 & -1 & 1 \\
$\dots$ & $\dots$ & $\dots$ \\
\end{tabular}
\end{center}

$a_n^0$ is either $1$ or $-1$ for any $n$.
Since $q_0$ is the root of polynomial $A_m$,
$$
\left( a_m^m q_0^{m-1} + \dots + a_m^1 \right) \cdot q_0 + a_m^0 = 0,
$$
and the expression in the brackets is a whole number.
Therefore, $q_0$ is the divisor of $a_m^0$. It is impossible if $q_0>1$. \qed

\vskip 1em
\noindent{\bf \fact{round}.\/}
$ \absolute{\nround{x} - x} \le 1/2$.

\noindent{\bf Proof.}
a) If $x \in \Z$, the left part is $0$. Let us consider cases $x \notin \Z$.

b) $\absolute{\nfloor{x} - x} < \absolute{\nceil{x} - x}$ and $\nround{x} = \nfloor{x}$
(see~\mref{round-sel}). Using~\mref{floorceil}, derive:
$$
\begin{array}{l}
x - \nfloor{x} < \nceil{x} - x = \nfloor{x} + 1 - x \quad\ergo \\
2 \absolute{x - \nround{x}} = 2x - 2\nfloor{x} < 1.
\end{array}
$$

c) $\absolute{\nfloor{x} - x} \ge \absolute{\nceil{x} - x}$ and $\nround{x} = \nceil{x}$.
Similar to case (b),
$$
x - (\nceil{x} - 1) \ge \nceil{x} - x \quad\ergo\quad
1 \ge 2 \absolute{\nround{x} - x}.
$$

In all 3 cases, the \fact{round} statement holds. \qed

\vskip 1em
\noindent{\bf \fact{atantaylor} Corollary.\/} $\forall\, x \in (0,1)$, $x - x^3/3 < \arctan x < x$.

\noindent{\bf Proof.\/}
Let $K=1$ and substitute $\eps$ with its estimation in~\mref{maclaurin2}.
Opening the inequality with absolute value,
$$
-\dfrac{x^5}{5} \le \arctan x - \left( x - \dfrac{x^3}{3} \right) \le \dfrac{x^5}{5},
$$
$$
\arctan x - x \le -\dfrac{x^3}{3} + \dfrac{x^5}{5}
    = -\dfrac{x^3}{3} \left( 1 - \dfrac{3x^2}{5} \right) < 0 \quad \forall\, x \in (0,1).
$$

Let $K=2$. Similar to the above,
$$
\arctan x - \left( x - \dfrac{x^3}{3} + \dfrac{x^5}{5} \right) \ge -\dfrac{x^7}{7},
$$
$$
\arctan x - \left( x - \dfrac{x^3}{3} \right) \ge
    \dfrac{x^5}{5} \left( 1 - \dfrac{5x^2}{7} \right) > 0
    \quad \forall\, x \in (0,1), \quad\hbox{\qed}
$$

\section{Conclusion}

We presented the method to build the simple Machin-like identity for any selected
starting term. We proved that Lehmer's measure of our identity can be made
arbitrarily small by selecting a single parameter. We presented the practical method, and
the Python program to compute parameters in our identity, along with the algorithm
for the actual computation of $\pi$ with pre-defined precision.

\section{Acknowledgements}

I am grateful to Dr.~Abrarov,~S.~M. for the discussion of the results.

\section{Special thanks}

\noindent
To my wife, Ekaterina, for the encouragement to put my findings on paper.

\noindent
To the open-source community for the tools to compute and create.

\end{document}